\documentclass[twoside,11pt]{article}
\usepackage{amsmath,amssymb, amsfonts, amstext,amsthm}

\newtheorem{thm}{Theorem}[section]
\newtheorem{den}[thm]{Definition}
\newtheorem{rek}[thm]{Remark}
\newtheorem{lea}[thm]{Lemma}
\newtheorem{asn}[thm]{Assumption}

\numberwithin{equation}{section}
\renewcommand{\newcommand}{\def}

\def\ba{\begin{array}}
\def\ea{\end{array}}
\def\be{\begin{eqnarray}}
\def\ee{\end{eqnarray}}
\def\bee{\begin{eqnarray*}}
\def\eee{\end{eqnarray*}}
\def\na{\nabla}

\def\al{\alpha}
\def\ep{\varepsilon}
\def\Ga{\Gamma}
\def\De{\Delta}
\def\na{\nabla}
\def\th{\theta}
\def\om{\omega}
\def\Om{\Omega}

\def\ph{\varphi}

\def\la{\lambda}
\def\fo{\forall}
\def\iy{\infty}

\def\si{\sigma}
\def\pp{\partial}
\def\de{\delta}

\begin{document}
\thispagestyle{plain}

\centerline{\LARGE\bf An impact of stochastic dynamic boundary}
\vskip0.2cm \centerline{\LARGE\bf conditions on the evolution of }
\vskip0.2cm \centerline{\LARGE\bf the   Cahn-Hilliard system
\footnote{\textbf{Date}: July 25, 2006. \\
\textbf{MSC 2000}: Primary-60H15, 37L55; secondary-35R60,37H10.\\
\textbf{Key words and phrases}: Stochastic partial differential
equations (SPDES), stochastic dynamic  boundary condition,
  nonlinear system under uncertainty, microscopic
mechanism on boundary, random dynamical system, random attractor.\\
\textbf{Acknowledgment}: We thank Alain Miranville and Dirk Blomker
for helpful comments. This work was partly supported by the NSF
Grants DMS-0209326 \& DMS-0542450. }}

\medskip

Abbreviated Title: \textbf{Impact of stochastic dynamic boundary
conditions}


\vskip1cm


\centerline{\large \bf Desheng Yang}

\centerline{\small\it
School of Mathematical Sciences and Computing
 Technology }
\centerline{\small\it Central South University, Changsha 410083,
China} \centerline{\small\it E-mail:dsyang@mail.csu.edu.cn}
\medskip
\centerline{and}
\medskip

\centerline{\large\bf Jinqiao Duan (\emph{Corresponding author}) }
\centerline{\small\it Department of Applied Mathematics, Illinois
Institute of Technology Chicago, IL 60616, USA}

\centerline{\small\it E-mail:duan@iit.edu}


\vskip1cm

{ \small

\noindent{\bf Abstract:}

Nonlinear systems are often subject to random influences.
Sometimes the noise enters the system through physical boundaries
and this leads to stochastic dynamic boundary conditions. A
dynamic, as opposed to static, boundary    condition  involves the
time derivative as well as spatial derivatives for the system
state variables on the boundary. Although stochastic \emph{static}
(Neumann or Dirichet type) boundary conditions have been applied
for stochastic partial differential equations, not much is known
about the dynamical impact of stochastic \emph{dynamic} boundary
conditions. The purpose of this paper is to study possible impacts
of stochastic dynamic boundary conditions on the long term
dynamics of the Cahn-Hilliard equation arising in the materials
science. We show that the dimension estimation of the random
attractor increases as the coefficient for the dynamic term  in
the stochastic dynamic boundary condition decreases. However, the
dimension of the random attractor is not affected by the
corresponding stochastic static boundary condition.




}

\section{Introduction}
\markboth{D. Yang and J. Duan }{ Impact of stochastic dynamic
boundary conditions }

The deterministic Cahn-Hilliard equation was introduced \cite{cah}
as a mathematical model for the description of phase separation
phenomena in materials such as binary alloys. The concentration
$\phi$ of one of the two components of the binary alloy satisfies
\be\label{1-1} \pp_t\phi =\De (-\De \phi+f(\phi)) \quad {\rm in
\;\; G},
\ee where domain $G:=\Pi_{i=1}^n(0, L_i), L_i>0, n\in\{1,2,3\}$ and
the boundary is denoted by $\Ga$. The function $f$
is the derivative of a logarithmic potential which is usually
approximated by a polynomial   with strictly positive dominant
coefficient.  Theoretical results on   asymptotic dynamical
behavior of   the system (\ref{1-1}), under either Neumann or
periodic boundary conditions, can be founded, for example, in
 the survey \cite{nov} or the book \cite{tem}, and \cite{dap1, blo, GMW}.

Recently physicists have considered phase separation phenomena in
confined systems \cite{Fischer1, Fischer2, Kenzler}, where
interactions with the walls have to be taken into account. This
leads to dynamic boundary conditions, i.e., $\pp_t \phi$ appears
in the boundary conditions. For example, one of the important
phase separation phenomena when the binary alloy is suddenly
cooled sufficiently is called spinodal decomposition. Once the
effective interaction between the wall (i.e., the boundary $\Ga$)
and two components of the alloy are short-ranged in spinodal
decomposition, the so-called dynamic boundary conditions have to
be taken into account, together with the the Cahn-Hilliard
equation. Mathematical results concerning the Cahn-Hilliard
equation with (deterministic) dynamic boundary conditions have
been obtained recently; see \cite{mir, rac, wuh}  and
  references therein.  We specially note that the so-called
viscous Cahn-Hilliard equation with \emph{dynamic} boundary
conditions
 \cite{mir},
\be\label{1-2}\ba{l}
\pp_t\phi=\De \mu, \; \mu=-\De\phi+\ep \pp_t\phi+f(\phi),\,{\rm in} \, G,\\
\pp_\nu \mu=0,\,{\rm on } \, {\Ga},\\
\pp_t\phi=\De_\|\phi-\la\phi-\pp_\nu\phi-g(\phi), \,{\rm on} \, {\Ga},\\
\phi(0)=\phi_0, \ea\ee where $\ep\ge 0$ is a small parameter,
$\De_\|$ is the Laplace-Beltrami operator on the boundary $\Ga$,
$\mu$ is chemical potential, $\nu$ is the unit outer normal vector
to $\Ga$ and $\la $ is some given positive constant. Miranville
and Zelik \cite{mir} have recently constructed a robust family of
exponential attractors for this system when the nonlinear
functions $f$ and $g$ are of arbitrary growths but satisfy some
dissipativity assumptions.

 Moreover, the environmental or
surrounding fluctuations may also influence the system evolution
and thus may be taken into account as well \cite{blo, dap1, EM}.
The present paper is concerned with the stochastic version of Eq.
(\ref{1-2}) with a stochastic dynamic boundary condition. It is
given by \be\label{1-4}\ba{l}
d\phi=\De \mu dt+ \si_1dW^{(1)},\; \mu=-\De\phi+\ep \pp_t\phi+f(\phi),\,{\rm in} \, G,\\
\pp_\nu \mu=0,\,{\rm on } \, {\Ga},\\
d\phi=(\De_\|\phi-\la\phi-\pp_\nu\phi-g(\phi))dt+\si_2dW^{(2)}, \,{\rm on} \, {\Ga},\\
\phi(0)=\phi_0, \ea\ee where $W^{(1)}$ and $W^{(2)}$ are
independent Wiener processes which will be explained in detail
later, and the constants $\si_i>0, i=1,2$, for the noise
intensities. The random fluctuation terms consisting of
$(W^{(1)},W^{(2)})$ act in the domain but also on the boundary
$\Ga$ and give a refined description of the underlying microscopic
mechanism in the phase separation phenomena. To simplify the
situation, we only consider the case with $g=0$, and the potential
$f(u)$ is a polynomial of odd degree with a positive leading
coefficient such as \be\label{1-3} f(u)=\sum_{k=1}^{2p-1} a_k u^k,
a_{2p-1}>0,p\ge2. \ee

There have only been a few works on spatially extended systems
under stochastic \emph{dynamic} boundary conditions
\cite{DuanGaoSchm, ChuSchm}. It is intriguing to know the impact
of such dynamic boundary conditions on the overall dynamics.


As in the deterministic case \cite{mir}, we introduce another
unknown function, namely, $\psi=\phi|_\Ga$, defined on the
boundary $\Ga$ and rewrite Eq.(\ref{1-4}) as a coupled system of
It\^o parabolic stochastic partial differential equations
(\textbf{SPDEs}) of the form: \be\label{1-5}\ba{l}
d\phi=\De \mu dt+\si_1 dW^{(1)}, \mu=-\De\phi+\ep \pp_t\phi+f(\phi),\,x\in G,\\
\pp_\nu \mu=0,\,{\rm on } \, {\Ga},\\
\phi(0)=\phi_0,\\
d\psi=(\De_\|\psi-\la\psi-\pp_\nu\phi)dt+\si_2 dW^{(2)}, \,  x\in \Ga ,\\
\psi(0)=\psi_0,\\
\phi|_\Ga=\psi.
\ea\ee

The aim of this paper is study the possible impact of the
stochastic dynamic boundary condition on the overall evolution of
the stochastic Cahn-Hilliard equation (\ref{1-5}). To this end, we
will, later on in Section 5 (see Remark 5.4), look at the dynamic
boundary condition but with a positive intensity parameter $\ep_0$
in front of the time derivative, \be \label{newbc} \frac{1}{\ep_0}
d\psi=(\De_\|\psi-\la\psi-\pp_\nu\phi)dt+\si_2 dW^{(2)}, \, x\in
\Ga. \ee We will be able to see an impact of the stochastic
dynamic boundary condition by looking at how the random attractor
varies with  the dynamic intensity parameter $\ep_0$.

Based the theory on   random dynamical systems  \cite{arn}, we
first define a cocycle by the stochastic flow corresponding to the
solution map for Eq.(\ref{1-5}). Due to the additive noise, by
introducing appropriate stochastic convolutions, we can solve the
problem (\ref{1-5}) pathwise, and thus construct a cocycle. As the
cocycle is given, we then prove the system possesses   a random
attractor. Since the system considered is non-autonomous, we
employ the pull-back technique to describe the asymptotic
behavior. Finally we estimate Hausdorff dimension  \cite{deb} of
the random attractor and obtain its dependence on the dynamic
intensity parameter $\ep_0$.

This paper is arranged as follows. In Section 2, we present some
preliminaries including some function spaces, inequalities, the
spectrum results on some operators, the definitions of the noises
and the basic framework of random dynamical systems. In Section 3,
we construct a cocycle, defined via the stochastic flow generated
by the stochastic Cahn-Hilliard equation (\ref{1-5}). In Section
4, we prove the existence of the random attractor. Finally, in
Section 5, we show that the dimension estimation of the random
attractor increases as the coefficient for the dynamic term  in
the stochastic dynamic boundary condition decreases. However, the
dimension of the random attractor is not affected by the
corresponding stochastic static boundary condition.

\section{Preliminaries}

In this section we introduce some function spaces, inequalities,
the known spectrum results on some operators \cite{tem, mir} and
the definitions of the noises \cite{dap2}. For convenience, we
also recall the basic framework of   random dynamical systems
theory \cite{arn} and random attractors.

\subsection {Function spaces and noises}

We denote by $H$ (respectively, $H(\Ga)$) the usual Sobolev space $L^2(G)$
(respectively, $L^2(\Ga)$) with the inner product $(\cdot, \cdot)$ (respectively,$(\cdot, \cdot)_\Ga$)
and the norm  $|\cdot|$(respectively, $(|\cdot|)_\Ga$). Let $m(u)$ be
the average over $G$
$$m(u)=\frac1{|G|}\int_G u(x) dx,$$
and $H_0$ be a subspace of $H$ defined by
$$H_0=\{u\in H: m(u)=0\}.$$

We define the linear unbounded operator $A=-\De$ with domain
$D(A)=\{u\in H^2(G)\cap H_0: \pp_\nu u=0\,{\rm on}\, \Ga \}$. Then
its spectrum consists of an infinite sequence of real eigenvalues
$0=\la_0<\la_1\le\la_2\le\dots\to \iy$. The corresponding
eigenvectors $\{e_i^{(1)} \} $ form a complete orthonormal basis
in $H$. Also notice that $\{e_i^{(1)}:i\in {\Bbb N} \}$ is a complete
orthonormal basis in $H_0$ and the function $e_0^{(1)}$ is equal to a
constant $|G|^{-1/2}$. By spectral theory, we may define the
fractional powers $A^s$ of $A$ , $s \in R$,  by
$$A^s u= \sum_{i=1}^\iy \la_i^s u_i e_i^{(1)} , $$
if $u=\sum_{i=0}^\iy  u_i e_i^{(1)}$. The domain of $A^{s/2}$ is
$$V_s=D(A^{s/2})=\{u=\sum_{i=0}^\iy  u_i e_i^{(1)}: \sum_{i=1}^\iy \la_i^s u_i^2<\iy\}.$$
This domain space is endowed with the seminorm $|u|_s=|A^{s/2}u|$,
the semiscalar product $(u, v)_s=(A^{s/2}u, A^{s/2}v)$ and the
norm:
$$\|u\|_s=(|u|^2_s+m^2(u))^{\frac12}.$$
Therefore, we have
$$V_0=H,A^0u=\bar u=u-m(u),\|\bar u\|_s=|\bar u|_s=|u|_s,$$
and the following Poincare type inequality:
$$|u|_s\le \la_1^{\frac{s-t}2}|u|_t, s\le t,u\in V_t,$$
together with the interpolation inequality:
$$ |u|_{\la s+(1-\la)t}\le |u|_s^\la |u|_t^{1-\la}, s\le t,\la \in [0,1],u\in V_t . $$

We also define the operator $B$ by the Lalace-Beltrami operator,
$B=-\De_\|$, with the domain $D(B)=H^2(\Ga)$. Its spectrum
consists of an infinite sequence of real eigenvalues
$0<\mu_1\le\mu_2\le\dots\to \iy$. The corresponding eigenvectors
$\{e^{(2)}_i \}_{i\in {\Bbb N}} $ form a complete orthonormal
basis in $H(\Ga)$. Like the Sobolev space $V_s$ defined on $G$, we
can define the space $V_{s,\Ga}$ on  $\Ga$ with the seminorm
$|\cdot|_{s,\Ga}$ and the semiscalar product $(\cdot,\cdot)_{s,
\Ga}$, and also have the corresponding interpolation inequality
and Poincare-type inequality.

For the additive stochastic term, we assume the following.

\begin{asn}\label{noise}
The stochastic process $W(t):=(W^{(1)}(t), W^{(2)}(t))$, defined
on a filtered probability space $(\Om,{\cal F},\{{\cal
F}_t\}_{t\in {\Bbb R}},P)$, is a two-side in time Wiener process
(Brownian motion) on $(H_0, H(\Ga))$, with covariance operator
$Q=(Q_1, Q_2)$. It is given as the expansions \be\label{2-1}
W^{(j)}(t)=\sum_{i=1}^\iy
\sqrt{\al_i^{(j)}}\beta_i^{(j)}(t)e^{(j)}_i,\quad {\rm with }\quad
Q_j e^{(j)}_i=\al_i^{(j)} e_i^{(j)}, j=1,2, \ee where
$\{(e_i^{(1)}, e_i^{(2)})\}_{i\in {\Bbb N}}$ are the orthonormal
basis of $(H_0, H(\Ga))$ from the eigenvectors of the operator $A$
and $B$. Moreover, $\{(\beta_i^{(1)}, \beta_i^{(2)})\}_{i\in {\Bbb
N} }$ are independent (two-sided in time) standard scalar Wiener
processes (Brownian motions)   on the probability space
$(\Om,{\cal F},P)$. It is also assumed that the covariance
operator $Q=(Q_1, Q_2)$ is a Hilbert-Schmidt operator, and is of
trace class, i.e. ${\rm tr}Q_j=\sum_{i=1}^\iy\al_i^{(j)}<\iy,
j=1,2$.
\end{asn}

The special case  $Q=I$, where $I$ is an identity operator, or
equivalently, $\al_i^{(j)}=1, j=1,2$ for all $i\in {\Bbb N}$,
corresponds to the case of a cylindrical white noise (white in
both time and space). Therefore, the cylindrical Wiener process
does not satisfy the above assumption and our result in this paper
does not apply to this case. In the deterministic case the
Cahn-Hilliard equation preserves the spatial average of $\phi$ in
time. According to the definition of $W(t)$, we have $\int_G
W^{(1)}(t,x) dx=0$. It guarantees that the white noise does not
destroy this property, i.e., the spatial average of $\phi$ is
still a conversed quantity: \be\label{2-2}
m(\phi(t))=m(\phi(0)).\ee

Let us consider the operator $A_\ep$ defined by
$$A_\ep u=(\ep+A^{-1})^{-1}Au,$$
for all $u\in D(A_\ep)=D(A)$. Then the spectrum of $A_\ep$
consists of eigenvalues $r_k$:
$$r_k=(\ep+\la_k^{-1})^{-1}\la_k=\frac{\la_k^2}{1+\ep\la_k},$$
and the corresponding eigenfunction $e_k^{(1)}$.

Under Assumption \ref{noise}, the linear equation
\bee\ba{l}
dz+\De^2z dt=\ep\De dz+\si_1 dW^{(1)}, \, x\in G,\\
\pp_\nu z|_\Ga=0,
\ea\eee
which can be rewritten as the abstract form:
\be\label{2-3}
(\ep+A^{-1})dz+Az dt=\si_1A^{-1} d W^{(1)}(t),\, {\rm in }\, H_0,
\ee
has a unique stationary solution given by
the It\^o integral and Fourier series expansion
\be\label{2-4}\ba{rl}
z(t)&=\si_1\int_{-\iy}^te^{-A_\ep(t-s)}(I+\ep A)^{-1} d W^{(1)}(s)\\
&=\si_1\sum_{k=1}^{+\iy}\frac{\sqrt{\al^{(1)}_k}}{\ep\la_k+1}\int_{-\iy}^t
e^{-r_k(t-s)} d\beta^{(1)}_k e^{(1)}_k, \ea\ee which is also
called as a stochastic convolution. As is well known, $z(t)$ is a
continuous Gaussian process in the space $H_0$. Moreover, we can
prove the process $z(t,x)$ satisfies the following regularity
(similar to \cite{dap1} and \cite{dap2}).

\medskip

Throughout the paper, the letters $C$ and $C_i$ denote some
generic constants, which may change from one line to another.

\begin{lea}\label{convolution}
Under Assumption \ref{noise}, the process  $\na z(t,x)$ has a version which is $\nu$-H\"older
continuous with respect to $(t,x)\in {\Bbb R}\times G$, for any $\nu\in [0,\frac14), \ep\ge0$.
Furthermore, if $\ep>0$, we have $Az(t,x)\in C^\nu({\Bbb R}\times G)$, for any $\nu\in[0,\frac12)$.
\end{lea}

\noindent{\it Proof.} For the sake of simplicity, we assume $G=(0,\pi)$.
Therefore, the eigenfunctions $\{e_i^{(1)}: i\in {\Bbb N}\}$ satisfy
$$|e_i^{(1)}|\le1, \;\;  |\na e_i^{(1)}|\le \sqrt{\la_i}.$$

From the representation (\ref{2-4}), we can write $\na z$ as
the following expansion
$$\na z(t,x)= \si_1\sum_{k=1}^{+\iy}\frac{\sqrt{\al^{(1)}_k}}{1+\ep\la_k}
\int_{-\iy}^t e^{-r_k(t-s)} d\beta^{(1)}_k \na e^{(1)}_k  $$
Therefore, using the inequality $|\na e^{(1)}_k(x)-\na e^{(1)}_k(y)|\le \la_k|x-y|$,
we have the estimate for any $x, y\in G, \ep\ge0$,

\be\label{2-5}\ba{rl}
E|\na z(t,x)-\na z(t,y)|^2
&\le \si_1 |x-y|^2\sum_{k=1}^{+\iy}\frac{\al^{(1)}_k\la_k^2}{(1+\ep\la_k)^2}
\int_{-\iy}^t e^{-2r_k(t-s)} d s  \\
&=\frac{\si_1}2|x-y|^2\sum_{k=1}^{+\iy}\frac{\al_k^{(1)}}{1+\ep\la_k}
\le \frac{\si_1}2 {\rm tr} Q_1|x-y|^2.
\ea\ee
Furthermore, using the inequality $|A e^{(1)}_k(x)-A e^{(1)}_k(y)|\le \la_k^{\frac32}|x-y|$,
for $\ep>0$, we obtain
\be\label{2-6}\ba{rl}
E|A z(t,x)-A z(t,y)|^2
&\le \si_1 |x-y|^2\sum_{k=1}^{+\iy}\frac{\al^{(1)}_k\la_k^3}{(1+\ep\la_k)^2}
\int_{-\iy}^t e^{-2r_k(t-s)} d s  \\
&=\frac{\si_1}2|x-y|^2\sum_{k=1}^{+\iy}\frac{\al_k^{(1)}\la_k}{1+\ep\la_k}
\le \frac{\si_1}{2\ep} {\rm tr} Q_1|x-y|^2.
\ea\ee

Moreover, fixing $t>s, x\in G$, we have
\bee\ba{l}
\quad E|\na z(t,x)-\na z(s,x)|^2\\
\le\si_1\sum_{k=1}^{+\iy}\frac{\al^{(1)}_k\la_k}{(1+\ep\la_k)^2} \{\int_s^te^{2r_k(\tau-t)} d
\tau+ \int_{-\iy}^s|e^{(\tau-t)r_k}-e^{(\tau-s)r_k}|^2d \tau\}\\
=\si_1\sum_{k=1}^{+\iy}\frac{\al^{(1)}_k\la_k}{(1+\ep\la_k)^2r_k}\{1-e^{r_k(s-t)}\}.
\ea\eee
Writing $C_\beta=\sup\limits_{x\ge0, y\ge0}\frac{e^{-x}-e^{-y}}{|x-y|^{2\beta}}$,
we get for $\beta\in [0,\frac14], \ep\ge0$,
\be\label{2-7}\ba{rl}
E|\na z(t,x)-\na z(s,x)|^2
&\le C\sum_{k=1}^{+\iy}\frac{\al^{(1)}_k\la_k}{(1+\ep\la_k)^2}r_k^{2\beta-1}|t-s|^{2\beta}\\
&=C\sum_{k=1}^{+\iy}\frac{\al^{(1)}_k\la_k^{4\beta-1}}{(1+\ep\la_k)^{2\beta+1}}|t-s|^{2\beta}\le C |t-s|^{2\beta}.
\ea\ee

Similarly, for $\ep>0, t>s, \beta \in [0,\frac12]$, we have
\be\label{2-8}
E|A z(t,x)-A z(s,x)|^2
\le C\sum_{k=1}^{+\iy}\frac{\al_k^{(1)}\la_k^{4\beta}}{(1+\ep\la_k)^{2\beta+1}}|t-s|^{2\beta}
\le C\ep^{-2\beta-1}|s-t|^{2\beta}.
\ee
Consequently, for $\beta \in [0,\frac14], \ep\ge0$, by (\ref{2-5})
and (\ref{2-7}), there exists a positive constant $C$ such that
\be\label{2-9}
E|\na z(t,x)-\na z(s,y)|^2\le C (|x-y|^2
+|t-s|^2)^\beta,
\ee
for $x, y\in G, t,s \in {\Bbb R}$.
Since the random variable $\na z(t,x)-\na z(s,y)$ is
Gaussian, then for all $m\in N$, we obtain
$$
E(|\na z(t,x)-\na z(s,y)|^{2m})\le C
(|x-y|^2+|t-s|^2)^{m\beta},
$$
which gives that $\na z(t,x)\in C^\nu({\Bbb R} \times G),$ for any $\nu\in [0, \frac14), \ep\ge0$
from Kolmogorov Theorem (\cite{oks}).
Similarly, by (\ref{2-6}) and (\ref{2-8}), for $\ep>0$, we have $Az(t,x)\in C^\nu({\Bbb R}\times G),
\nu\in [0,\frac12)$.
$\blacksquare$

\begin{rek}
Note that the estimates for eigenfunctions, $|e_i^{(1)}|\le 1$ and
$|\na e_i^{(1)}|\le \sqrt{\la_i}$,  are  true for rectangular and
triangular domains, but not true for disks. It is unknown which
geometrical conditions (on the domain) imply these estimates. See
detailed discussions in \cite{Desi, Blomker}.
\end{rek}

\subsection{Random dynamical systems}

Let $(\Om,{\cal F},P)$ be a probability space and $\{\th_t:\Om\to
\Om,t\in {\Bbb R}\}$ a family of measure-preserving transformations such
that $\th_0=I$, and $\th_{t+s}=\th_t\circ \th_s$, for all $t,s\in
R$. Here $\{\th_t\}$ is called   a metric or \emph{driving}
dynamical system on $(\Om,{\cal F},P)$. We always assume that
$\th$ is ergodic under probability measure $P$.

\begin{den}\label{rds}  Let $(X,d)$ be a Polish space
(i.e., complete separable metric space).  A measurable map
$$\ph: {\Bbb R}\times\Om\times X\to X,\;\; (t,\om,x)  \mapsto \ph  (t,\om)x, $$
is called a random dynamical system (\textbf{RDS}) if  $\ph $
satisfies the cocycle property: $\ph (0,\om)=I,\; \ph
(t+s,\om)=\ph (t,\th_s\om)\ph (s,\om)$, for all $t, s\in {\Bbb R}$ and
P-a.s. $\om\in\Om$.
\end{den}

A RDS is said to be continuous if $\ph (t,\om):X\to X$ is
continuous, P-a.s. for $\omega$ and for every $t\in {\Bbb R}$. Notice
that the RDS $\ph$ defined above is two-sided in time and
invertible, $\ph(t,\om)^{-1}=\ph(-t,\th_t\om)$ P-a.s. As random
attractors are characterized by random sets, we have to deal with
some new concepts such as absorbing sets, attraction and
invariance. Specially the pull-back approach, starting from $-\iy$
and observing time $0$, has been extensively explored in the
theory of random dynamical systems. Before defining them, we
recall   a compact random set as defined below.

\begin{den} A set-valued map $K:\Om\to 2^X$, the set of all subsets of  $X$, is called a random compact set if
$K(\om)$ is a compact P-almost surely and if $\om\to d(x,K(\om))$
is measurable for each $x\in X$, where $d(x,M):=\inf_{y\in
M}d(x,y)$.
\end{den}

The pull-back technique is used to define the following concepts.

\begin{den}\label{random attractor} Let $A(\om)$ and $B(\om)$ be two random sets. We say

{\narrower
\item{(1).}  $A(\om)$ attracts $B(\om)$ if P-a.s.,
$$\lim_{t\to \iy } {\rm dist }( \ph (t,\th_{-t}\om)B(\th_{-t}\om),A(\om))=0, $$
where ${\rm dist}(\cdot,\cdot)$ denotes the Hausdorff semidistance in $X$, defined as
${\rm dist} (A,B)=\sup_{x\in A}\inf_{y\in B} d(x,y), for \fo A,B\subset X$.
\item{(2).} $A(\om)$ absorbs $B(\om)$ if there exists $t_B(\om)$ such that
for all $t\ge t_B(\om)$,
$$\ph(t,\th_{-t}\om)B(\th_{-t}\om)\subset A(\om) . $$
holds P-a.s..

}
\end{den}

It is clear that a random absorbing set is attracting. We call a
invariant compact attracting set as a random attractor. Formally
 \cite{cra1, sch1},

\begin{den}\label{att} A random set ${\cal A}(\om)$ is called
 a random attractor for the RDS $\ph $ if
P-a.s.

{\narrower
\item{(1).}  ${\cal A}(\om)$ is a random compact set.
\item{(2).}  ${\cal A}(\om)$ is strictly invariant, i.e., $\ph (t,\om){\cal A}(\om)={\cal A}(\th_t\om)$,
for $\fo t\ge0$.
\item{(3).}  ${\cal A}(\om)$ attracts all bounded deterministic sets in $X$.

}
\end{den}

The existence result of random attractors is stated as follows
 (e.g.,\cite{cra1, cra3, sch1, sch2}).

\begin{thm}\label{theorem} If there exists a random compact set absorbing every
bounded non-random set $B\subset X$, then the RDS $\ph$ possesses
a random attractor ${\cal A}(\om)$
$${\cal A}(\om)=\overline{\underset {B\subset X}  \bigcup \Lambda_B(\om)}, $$
where
$\Lambda_B(\om):=\underset{s\ge0}\bigcap\overline{\underset{t\ge
s}\bigcup\ph (t,\th_{-t}\om)B}$ is the omega-limit set of $B$.
\end{thm}

\begin{rek}\rm
\item{(1).} In fact, a random attractor $\cal A(\om)$ in
Definition \ref{random attractor} is a global random set attractor
which is uniquely determined by a attracting compact set
(\cite{cra2}). \item{(2).} $\ph (t,\th_{-t}\om)x$ can be
interpreted as the position at $t=0$ of the trajectory which was
at $x$ at time $-t$, that is, while time $t$ is moving, the
trajectory $\ph (t,\th_{-t}\om)x$ is always at the position at
time zero. Therefore, the random attractor in Definition \ref{att}
is also called as the ``pullback attractor".
\end{rek}

One of the results in the theory of global attractors for
deterministic systems is that the Hausdorff dimension of the
attractor is often finite.
Although the random attractor is not uniformly bounded, it is
expected that the techniques on the Hausdorff dimension of a
global attractor of a deterministic system can be generalized to
the stochastic case under some assumptions \cite{deb, sch2}. In
fact, based on Laypunov exponents,  Debussche \cite{deb} showed
that the Hausdorff dimension of a random attractor is finite if
the corresponding cocycle $\ph(t,\om)$ satisfies some properties,
especially uniformly differentiability. The following conclusion
is due to Debussche  \cite{deb}.

\begin{thm}\label{dimension} Let ${\cal A(\om)}$ be a compact measurable
set which is invariant under a random map
$S(\om),\om\in \Om$, for some ergodic metric dynamical system
$(\Om,{\cal F},P,(\th_t)_{t\in {\Bbb R}})$. Assume that
\begin{enumerate}
\item   $S(\om)$ is almost surely uniformly differentiable on ${\cal
A(\om)}$,that is, for every $u,u+h \in {\cal A(\om)}$, there
exists $D(S(\om;u))$ in $ {\cal L}(H) $,  the space of bounded
linear operator from a Hilbert space $H$ to itself, such that
$$|S(\om)(u+h)-S(\om)u-DS(\om;u)h|\le \bar k(\om)|h|^{1+\mu},$$
where $\mu>0, \bar k(\om)$ is a random variable satisfying $\bar
k(\om)\ge1, E(\log \bar k)<\iy$.
\item $\om_d(DS(\om,u))\le  \bar\om_d(\om)$ for $u\in {\cal A(\om)}$ and
some random variable $\bar\om_d(\om)$ satisfying
$E(\log(\bar\om_d))<0$, where
$$\om_d(L)=\al_1(L)\cdots\al_d(L),\al_i(L)=\inf_{{F\subset H}\atop{{\rm dim} F\le i-1}}\sup_{{\ph\in F^\perp}
\atop{ |\ph|_H=1}}|L\ph|_H\quad {\text for}\quad L\in {\cal L}(H).$$
\item $\al_1(DS(\om,u))\le \bar \al_1(\om)$, for $u\in{\cal A}(\om)$ and
a random variable $\bar \al_1(\om)\ge 1$ with
$E(\log\bar\al_1)<\iy$.
\end{enumerate}
Then the Hausdorff dimension $d_H({\cal A(\om)})$ of ${\cal
A(\om)}$ is less than $d$ almost surely.
\end{thm}

\section{Stochastic flow}

The stochastic Cahn-Hilliard equation with stochastic dynamical
boundary conditions  is non-autonomous, and  thus it is impossible
to define a semigroup on the phase space. In fact, a solution of
the stochastic equation gives a stochastic flow instead of a
semigroup. In this section, we introduce stochastic convolutions
to solve the problem (\ref{1-5}) pathwise and obtain the
corresponding stochastic flow. In fact, the stochastic flow
satisfies the so-called cocycle property and thus leads to a RDS
modelling our system.

Let $\Om$ be the set of continuous functions
with value $(0,0)\in {\Bbb R}^2$ at 0
$$\Om=\{\om\in C({\Bbb R},{\Bbb R}^2):\om(0)=(0,0)\}.$$
Let ${\cal F}$ be the Borel sigma-algebra induced by the
compact-open topology of $\Om$, and  $P$ a Wiener measure on
$(\Om,{\cal F})$. We write
$$W(t)=(W^{(1)}(t),W^{(2)}(t))=(\om^{(1)}(t), \om^{(2)}(t))=\om(t), t\in {\Bbb R}, \om\in \Om,$$
and define
\bee
\th_t\om(s)=\om(t+s)-\om(t),\quad t\in {\Bbb R}.
\eee
It is easy to get $\th_t\circ\th_s=\th_{t+s}$ in terms of (\ref{noise}).
Thus $(\Om,{\cal F},P,(\th_t)_{t\in {\Bbb R}})$
is an ergodic metric dynamical system which models the white
noise.

Following \cite{mir}, we introduce the weak energy space $\Bbb
L_\ep$ with the norm
$$|(\phi, \psi)|_{\Bbb L_\ep}^2=\ep |\phi|^2+\|\phi\|^2_{-1}+|\psi|^2_\Ga,$$
and the phase space $\Bbb V$ with the norm
$$|(\phi, \psi)|_{\Bbb V}^2=\|\phi\|^2_1+\|\psi\|^2_{1,\Ga}.$$

For convenience, using the operator $A^{-1}:H_0\to H_0$, we
rewrite the system (\ref{1-5}) as follows: \be\label{3-1}\ba{l}
(\ep+A^{-1})d\phi=(\De \phi-<\pp_\nu\phi>_\Ga-\bar f(\phi)) dt+\si_1A^{-1} dW^{(1)},\,x\in G,\\
\phi(0)=\phi_0,\\
d\psi=(\De_\|\psi-\la\psi-\pp_\nu\phi)dt+\si_2 dW^{(2)}, \,  x\in \Ga ,\\
\psi(0)=\psi_0,\\
\phi|_\Ga=\psi,
\ea\ee
where $<u>_\Ga=\frac1{|G|}\int_\Ga u(x) dx$ and $\bar f=f-m(f)$.

Let $z^{(1)}$ and $z^{(2)}$ be stationary solutions of (\ref{2-3}) and the linear equation
$$d z=(\De_\|z-\la z) dt+\si_2 dW^{(2)},\, x\in \Ga,$$
respectively. Similar to Lemma \ref{convolution}, we can obtain
$\na_\| z^{(2)}(t, x)\in C^\nu({\Bbb R}\times \Ga), \nu\in [0,\frac14)$.

Set $(u, v)=(\phi-z^{(1)}, \psi-z^{(2)}).$ Then $(u,v)$ satisfies
\be\label{3-2}\ba{l}
(\ep+A^{-1}) \pp_t u=\De u-<\pp_\nu u>_\Ga-\bar f(u+z^{(1)}),\, x\in G,\\
u(0)=\phi_0-z^{(1)}(0),\\
\pp_t v=\De_\| v-\la v-\pp_\nu u, \,  x\in \Ga ,\\
v(0)=\psi_0-z^{(2)}(0),\\
u|_\Ga=v+z^{(2)}-z^{(1)}|_\Ga,
\ea\ee
where we have used the fact $\pp_\nu z^{(1)}|_\Ga=0$.

According to \cite{mir}, the system (\ref{3-2}) possesses a unique solution
$(u,v)\in C([s,t], {\Bbb L_\ep})$, for
every initial value $(u(s), v(s))\in {\Bbb L_\ep}, s\in {\Bbb R}, t\ge s, \om \in \Om$, and so, there exists
a continuous operator $S_\ep (t,s;\om)$ on the weak energy space $\Bbb L_\ep$:
$$S_\ep(t,s;\om): {\Bbb L_\ep}\to {\Bbb L_\ep}, S_\ep(t,s;\om)(\phi(s),\psi(s))=(\phi(t), \psi(t)).
$$
The corresponding stochastic
flow can be defined by
$$\ph(t,\om)(\phi(0),\psi(0))=S_\ep(t,0;\om)(\phi(0),\psi(0)).$$
Notice that P-a.s.,
$$S_\ep(t,s;\om)=S_\ep(t,r;\om)S_\ep(r,s;\om), S_\ep(t,s;\om)=S_\ep(t-s,0;\th_s\om).$$
It implies that $\ph $ satisfies the cocycle property. In the end,
$\ph$ gives a continuous random dynamical system on $\Bbb L_\ep$
over $(\Om,{\cal F},P,(\th_t)_{t\in {\Bbb R}})$  associated with the
stochastic Cahn-Hilliard equation (\ref{3-1}) under stochastic
dynamic boundary condition.

\section{The random attractor }

In this section, we prove the stochastic Cahn-Hilliard equations
(\ref{3-1}) possesses a random attractor. We first recall that the
stochastic system (\ref{3-1}) possesses the conservation law
(\ref{2-2}) and, consequently, we cannot expect to construct a
random attractor in the whole phase space $\Bbb L_\ep$. Therefore,
we work in the affine space
$$ \Bbb L_\ep^\beta=\{u\in \Bbb L_\ep: m(u)=\beta, 0<\ep\le1 \},$$
with the norm $\Bbb L_\ep$, and prove the restriction of the
cocycle $(\phi, \psi)$ on the $\Bbb L_\ep^\beta$ possesses a
compact absorbing set. By Theorem  \ref{theorem}, we conclude that
the stochastic Cahn-Hilliard equations (\ref{3-1}) with stochastic
dynamic boundary conditions has a random attractor.

We begin by proving the existence of a absorbing set in  $\Bbb
L_\ep^\beta$.

\begin{lea}\label{absorb} Given any ball of $\Bbb L_\ep$, $B(0,\rho)$ centered at $0$ of radius $\rho$, for
any $-1\le t\le 1, \om \in \Om$, there exist random variables $R_t(\om)$ and
$t(\rho;\om)<-2$ such that for any $s\le
t(\rho;\om),(\phi_s, \psi_s)\in B(0,\rho), $
\be\label{4-1} |S_\ep(t,s;\om)(\phi_s, \psi_s) | \le R_t(\om), \ee
holds P-a.s.
\end{lea}
\noindent {\it Proof.} Let $(u(t),v(t))=(u(t,\om; s,u_s), v(t,\om;
s,v_s))$ be the solution to Eqs.(\ref{3-2})  with the initial
value $(u_s, v_s)=(\phi_s-z^{(1)}(s),\psi_s-z^{(2)}(s))$.
Multiplying the first equation of (\ref{3-2}) by $\bar u$ and
integrating over $G$, we obtain \be\label{4-2} \frac12
\pp_t(\ep|\bar u|^2+|u|_{-1}^2)=(\De u,\bar u)-(\bar f(u+z^{(1)}),
\bar u), \ee where we have used $m(\bar u)=0$. We denote $\bar
v=v-m(u)$ and estimate the first term of the right hand side of
Eq.(\ref{4-2}) as follows: \bee\ba{rl}
(\De u,\bar u)&=(-\pp_t v+\De_\|v-\la v, \bar v)_\Ga+(\pp_\nu u, z^{(2)}-z^{(1)})_\Ga-|u|_1^2\\
&=-|u|_1^2-\frac 12\pp_t|\bar v|_\Ga^2-|\na_\| v|_\Ga^2-\la |v|^2_\Ga+\la(v, m(u))_\Ga+(\pp_\nu u,z^{(2)}-z^{(1)})_\Ga.
\ea\eee
Notice that
$$\la(v, m(v))_\Ga\le \frac\la2|v|_\Ga^2+C\beta^2,$$
and
$$\ba{rl}
(\pp_\nu u, z^{(2)}-z^{(1)})_\Ga
&\le \|\pp_\nu u\|_{-\frac12, \Ga}\|z^{(2)}-z^{(1)}\|_{\frac12, \Ga}\\
&\le
\frac14\|u\|_1^2+(\|z^{(2)}\|_{\frac12,\Ga}+\|z^{(1)}\|_{\frac32})^2.
\ea$$ Therefore, we have \be\label{4-3}\ba{l} (\De u,\bar u)\le
-\frac 12\pp_t|\bar v|_\Ga^2-|\na_\| v|_\Ga^2-\frac\la2 |v|^2_\Ga
-\frac34|u|_1^2+C\beta^2+(\|z^{(2)}\|_{\frac12,\Ga}+\|z^{(1)}\|_{\frac32})^2.
\ea\ee
Since $a_{2p-1}>0$, we see that
$$f(x)x\ge\frac12a_{2p-1}x^{2p}-C, |f(x)|\le 2 a_{2p-1}|x|^{2p-1}+C,$$
for some positive constant $C$. Therefore, we have the estimate on
the second term of the right hand side of Eq.(\ref{4-2}): \bee
\ba{rl} -(\bar f(u+z^{(1)}), \bar u)
&=-(f(u+z^{(1)}), u+z^{(1)}-m(u)-z^{(1)})\\
&\le -\frac{a_{2p-1}}2|u+z^{(1)})|_{L^{2p}}^{2p}+(f(u+z^{(1)}),m(u)+z^{(1)})+C,
\ea\eee
and
$$\ba{rl}
(f(u+z^{(1)}),m(u)+z^{(1)})
&\le 2a_{2p-1}\int_G(|u+z^{(1)}|^{2p-1}+C)|z^{(1)}+m(u)| dx\\
&\le \frac18a_{2p-1}|u+z^{(1)}|_{L^{2p}}^{2p}+C(|z^{(1)}|^{2p}_{L^{2p}}+\beta^{2p}+1).
\ea$$
In the end, we obtain
\be\label{4-4}
-(\bar f(u+z^{(1)}), \bar u)\le -\frac{3a_{2p-1}}8 |u+z^{(1)}|_{L^{2p}}^{2p}
+C(|z^{(1)}|^{2p}_{L^{2p}}+\beta^{2p}+1).
\ee
Summing (\ref{4-3}) and (\ref{4-4}) into (\ref{4-2}), for $\ep\in (0,1]$, we obtain
\be\label{4-5}\ba{l}
\quad \pp_t(\ep|\bar u|^2+|u|_{-1}^2+|\bar v|_\Ga^2)+\frac {\la_1\ep}4|\bar u|^2+\frac{\la_1^2}4|u|_{-1}^2+\la|v|_\Ga^2\\
\quad +|u|_1^2+2|v|_{1,\Ga}^2+\frac{3a_{2p-1}}4|u+z^{(1)}|^{2p}_{L^{2p}}\\
\le C(\|z^{(1)}\|_{\frac32}^{2p}+\|z^{(2)}\|^2_{\frac12, \Ga}+\beta^{2p}+1)
\ea\ee
where we have used inequalities $|u|_{-1}\le\la_1^{-1}|u|_1$,
$|\bar u|^2\le \la_1^{-1}|u|_1^2$ and
the Sobolev embedding theorem $V_{\frac32}\hookrightarrow L^{2p} (n=1,2,3)$.

Denoting $\kappa=\min\{\frac{\la_1}4, \frac{\la_1^2}4, \la\}>0$ and applying the Grownwall lemma to
(\ref{4-5}), we obtain
\be\label{4-6}\ba{l}
\quad |(\bar u(t), \bar v(t))|_{\Bbb L_\ep^\beta}^2+
\int_s^te^{-\kappa(t-\tau)}(|(u(\tau), v(\tau))|_{\Bbb V}^2+\frac{3}4a_{2p-1}|u+z^{(1)}|^{2p}_{L^{2p}}) d \tau\\
\le|(\bar u(s), \bar v(s))|_{\Bbb L_\ep^\beta}^2e^{-\kappa(t-s)}
+C\int_s^te^{-\kappa(t-\tau)}(\|z^{(1)}\|_{\frac32}^{2p}+\|z^{(2)}\|^2_{\frac12, \Ga}+\beta^{2p}+1) d \tau.
\ea\ee
Thus, for any ball $B(0,\rho)\subset {\Bbb L}_\ep^\beta$ and $(\phi(s), \psi(s))\in B(0,\rho)$, there
exists a random variable $t(\rho;\om)<-2$ such that for all $s<t(\rho;\om), t\in [-1,1]$,
\be\label{4-7}\ba{l}
\quad |(u(t),v(t))|_{\Bbb L_\ep^\beta}^2+\int_s^t|(u(\tau), v(\tau))|_{\Bbb V}^2 d \tau
+\int_s^t|u(\tau)+z^{(1)}(\tau)|^{2p}_{L^{2p}}d\tau\\
\le C(1+\beta^{2p}+|(z^{(1)}(t),z^{(2)}(t))|^2_{\Bbb L_\ep}
+\int_{-\iy}^te^{-\kappa(t-\tau)}(\|z^{(1)}\|_{\frac32}^{2p}+\|z^{(2)}\|^2_{\frac12,
\Ga}) d \tau). \ea\ee Denote by $r_t^2(\om)$ the right hand side
of the above inequality, which is finite P-a.s. by Lemma
\ref{convolution}. Writing
$R_t^2(\om)=2(r_t^2(\om)+|(z^{(1)}(t),z^{(2)}(t))|^2_{\Bbb
L_\ep})$, we complete  the proof of the lemma. $\blacksquare$

\medskip

Lemma \ref{absorb} shows that for any deterministic bounded set $B\subset B(0,\rho)$ in $\Bbb L_\ep^\beta$,
there exists a random time $t(\rho;\om)<-2$ such that for any $s<t(\rho;\om)$,
$$S_\ep(-1,s;\om)B\subset B(0,R_{-1}(\om)).$$
Noticing that
$$\ph(-s,\th_s\om)=S_\ep(-s,0;\th_s\om)=S_\ep(0,s;\om)=S_\ep(0,-1;\om)S_\ep(-1,s;\om),$$
we find that ${\cal B}(\om):=S_\ep(0,-1;\om)B(0,R_{-1}(\om))$ is a random absorbing set in ${\Bbb L_\ep^\beta}$.
Moreover, we can prove ${\cal B}(\om)$ is compact. In order to do it, we first give some estimates.

Multiplying the first equation of (\ref{3-2}) by $A u$ and
integrating over $G$, we obtain \be\label{4-8}
\frac12\pp_t(\ep|u|_1^2+|\bar u|^2)+|u|_2^2=-(\bar f(u+z^{(1)}),
Au). \ee
It is possible to find a positive constant $C$ such that
$$ f'(u)>\frac{2p-1}2a_{2p-1}u^{2p-2}-C,$$
and
$$
|f'(u)|\le 2(2p-1)a_{2p-1}u^{2p-2}+C,
 $$
and thus, it gives
$$\ba{l}
\quad -(f(u+z^{(1)}), Au)\\
=-\int_G \na f(u+z^{(1)}) \na u dx \\
\le -\frac{1}4(2p-1)a_{2p-1}\int_G|u|^{2p-2}|\na u|^2 dx+C|u|_1^2+
C_1(|u+z^{(1)}|_{L^{2p}}^{2p}+|\na z^{(1)}|_{L^{2p}}^{2p}).
\ea$$
Finally, (\ref{4-8}) yields
\be\label{4-9}\ba{l}
\quad \pp_t(\ep|u|_1^2+|\bar u|^2)+2|u|_2^2+\frac{(2p-1)a_{2p-1}}2\int_G|u|^{2p-2}|\na u|^2 dx\\
\le C|u|_1^2+
C_1(|u+z^{(1)}|_{L^{2p}}^{2p}+|\na z^{(1)}|_{L^{2p}}^{2p}).
\ea\ee
For any $t\in [-1,0], s<-1$, we have
$$\ba{l}
\quad \ep|u(t)|_1^2+|\bar u(t)|^2\\
\le\ep|u(s)|_1^2+|\bar u(s)|^2+C\int_s^t |u(\tau)|_1^2 d\tau
+C_1\int_s^t(|u(\tau)+z^{(1)}(\tau)|_{L^{2p}}^{2p}+|\na z^{(1)}(\tau)|_{L^{2p}}^{2p})d\tau.
\ea$$
Integrating over $[-2,-1]$ on $s$ leads to
\be\label{4-10}
\ep|u(t)|_1^2+|\bar u(t)|^2
\le C\int_{-2}^{-1} |u(t)|_1^2 dt
+C_1\int_{-2}^{-1}(|u+z^{(1)}|_{L^{2p}}^{2p}+|\na z^{(1)}|_{L^{2p}}^{2p})dt,
\ee
which yields
$$\ba{rl}
\int_{-1}^0|u(t)|_2^2 dt&\le \ep|u(-1)|_1^2+|\bar u(-1)|^2+C\int_{-1}^0 |u(t)|_1^2 dt
+C_1\int_{-1}^0(|u+z^{(1)}|_{L^{2p}}^{2p}+|\na z^{(1)}|_{L^{2p}}^{2p})dt\\
&\le C\int_{-2}^{0} |u(t)|_1^2 dt
+C_1\int_{-2}^{0}(|u+z^{(1)}|_{L^{2p}}^{2p}+|\na z^{(1)}|_{L^{2p}}^{2p})dt.
\ea
$$
Therefore, by (\ref{4-7}), for $s<t(\rho;\om)$, we obtain
$$\sup_{t\in[-1,0]}|u(t,\om;s,u_s)|<\iy, \int_{-1}^0|u(t,\om;s,u_s)|_2^2 dt<\iy.$$
The Sobolev embedding theorem $V_2\hookrightarrow L^r$ for $r\ge
2, n<4$, and the interpolation approach imply that for $q\in
[2,\iy)$ and $\beta\in [2,\iy)$, \be\label{4-11}
\int_{-1}^0|u(t)|_{L^q}^\beta dt<\iy, {\rm and\,\, thus},
\int_{-1}^0|\phi (t)|_{L^q}^\beta dt<\iy. \ee
Now we state that
${\cal B}(\om)$ satisfies the following compact property.

\begin{lea}\label{compact}
The random set ${\cal B}(\om)$ described above is a compact absorbing set, that is, it is compact
and absorbs any non-random bounded set:
for every bounded deterministic set $B\subset
B(0,\rho)$, we have
$$\ph(-s,\th_s\om)B\subset {\cal B}(\om) ,$$
holds P-a.s., for any $s\le t(\rho;\om)$, $t(\rho;\om)$ defined in
Lemma \ref{absorb}.
\end{lea}
\noindent{\it Proof.}  We only need to prove that the random absorbing set
${\cal B}(\om)$ is compact. Let $\{(\phi_0^n, \psi_0^n):n\in {\Bbb N}\}$ be a sequence in ${\cal B}(\om)$
and $(u^n, v^n)$ a solution to Eqs.(\ref{3-2}) with the initial value
$(u^n(0), v^n(0))=(\phi_0^n-z^{(1)}(0), \psi_0^n-z^{(2)}(0))$. We easily get, from (\ref{4-5}), that
for any $t\in [-1,0]$
\bee\ba{l}
\quad |(u^n(t),v^n(t))|_{\Bbb L_\ep^\beta}^2+
\int_{-1}^t|(u(\tau), v(\tau))|_{\Bbb V}^2 d \tau\\
\le|(u(-1), v(-1))|_{\Bbb L_\ep^\beta}^2
+C(1+\beta^{2p}+\int_{-1}^t (\|z^{(1)}\|_{\frac32}^{2p}+\|z^{(2)}\|^2_{\frac12, \Ga}) d \tau).
\ea\eee
Since $|(u(-1), v(-1))|_{\Bbb L_\ep^\beta}^2\le r^2_{-1}(\om)$, we deduce that $\{(u^n, v^n):n\in {\Bbb N}\}$
is bounded in
$L^\iy(-1,0; {\Bbb L_\ep^\beta})\cap L^2(-1,0; {\Bbb V})$, and so, it is compact in
$L^2(-1,0; {\Bbb L_\ep^\beta})$. Hence, there exists a
subsequence $\{(u^{n_k}, v^{n_k}) :n\in {\Bbb N}\}$
convergent to a function $(u,v)$ in $L^2(-1,0;{\Bbb L_\ep^\beta})$.
Moveover $(u,v)$ is a solution of Eqs. (\ref{3-2}). Let
$(\phi_0, \psi_0)=(u(0)+z^{(1)}(0), v(0)+z^{(2)}(0))$, it is easy to yield
$$\phi_0-\phi_0^{n_k}=u(0)-u^{n_k}(0), \psi_0-\psi_0^{n_k}=v(0)-v^{n_k}(0). $$
We now prove the subsequence $\{(\phi_0^{n_k}, \psi^{n_k}_0)\}$ is convergent to $(\phi_0, \psi_0)$ in
${\Bbb L_\ep^\beta}$. Writing
$(X(t), Y(t))=(u^{n_k}(t)-u(t), v^{n_k}(t)-v(t))$, then $(X(t), Y(t))$ satisfies the following equations
\be\label{4-12} \ba{l}
(\ep+A^{-1}) \pp_t X=\De X-<\pp_\nu X>_\Ga-(\bar f(u^{n_k}+z^{(1)})-\bar f(u+z^{(1)})),\, x\in G,\\
X(0)=\phi_0^{n_k}-\phi_0,\\
\pp_t Y=\De_\| Y-\la Y-\pp_\nu X, \,  x\in \Ga ,\\
Y(0)=\psi_0^{n_k}-\psi_0,\\
X|_\Ga=Y. \ea\ee Multiplying the first equation of (\ref{4-12}) by
$\bar X$ and integrating over $G$, we obtain \be\label{4-13}
\frac12\pp_t(\ep|\bar X|_1^2+|X|_{-1}^2)=(\De X,\bar X)-(\bar
f(u^{n_k}+z^{(1)})-\bar f(u+z^{(1)}),\bar X). \ee The two terms of
the right hand side of (\ref{4-13}) satisfy the following
estimates, respectively,
$$\ba{l}
(\De X, \bar X)\le -\frac12\pp_t|\bar Y|_\Ga^2-|\na_\|Y|_\Ga^2-\frac\la2|Y|_\Ga^2-|X|_1^2+Cm^2(X),\\
-(\bar f(u^{n_k}+z^{(1)})-\bar f(u+z^{(1)}),\bar X)\\
\le C|X|^2+C_1|m(X)||X|(1+|u^{n_k}+z^{(1)}|_{4p-4}^{4p-4}+|u+z^{(1)}|_{4p-4}^{4p-4}),
\ea$$
where we have used the fact $f'(u)\ge -C$, and, consequently, from (\ref{4-12}) and the
estimate $|m(X)|\le C|X|$, we deduce
\be\label{4-14}\ba{rl}
\pp_t|(X(t), Y(t))|_{\Bbb L_\ep^\beta}^2
&\le C(1+|u^{n_k}+z^{(1)}|_{L^{4p-4}}^{4p-4}+|u+z^{(1)}|_{L^{4p-4}}^{4p-4})|X|^2\\
&\le C(1+|u^{n_k}+z^{(1)}|_{L^{4p-4}}^{4p-4}+|u+z^{(1)}|_{L^{4p-4}}^{4p-4})|X|_{\Bbb L_\ep^\beta}^2\\
&:=C(t)|X|_{\Bbb L_\ep^\beta}^2,
\ea\ee
which leads to
$$|(X(0), Y(0))|_{\Bbb L_\ep^\beta}^2\le |(X(t), Y(t))|_{\Bbb L_\ep^\beta}^2 e^{\int_t^0 C(s) ds}.
$$
Integrating the above inequality over $[-1,0]$, we obtain
\be\label{4-15}
|(X(0), Y(0))|_{\Bbb L_\ep^\beta}^2\le e^{\int_{-1}^0 C(s) ds}|(X(t), Y(t))|^2_{L^2(-1,0; \Bbb L_\ep^\beta)},
\ee
which implies that the sequence  $\{(\phi_0^{n_k}, \psi^{n_k}_0)\}$ converges to $(\phi_0, \psi_0)$, and thus,
${\cal B}(\om)$ is compact. $\blacksquare$

Finally, using Theorem \ref{theorem}, we conclude
\begin{thm}
The random dynamical system associated with the stochastic
Cahn-Hilliard equation (\ref{1-5}) with stochastic dynamic
boundary conditions possesses a random attractor ${\cal A_\ep
(\om)}$.
\end{thm}

A compact absorbing set with the $\om$-wise attraction property guarantees
the existence of the random attractor, but the union in $\om$
of ${\cal A_\ep (\om)}$ is not compact in general.
Since $P$ is invariant under $\th_t$, the forward attraction property from
0 to $\iy$ can be obtained w.r.t. convergence in probability,
that is,
$$\lim_{t\to\iy}P({\rm dist}(\ph(t,\om)B,{\cal A_\ep }(\th_t\om))<\eta)=1,$$
holds for any $\eta>0$ and all deterministic bounded set $B\subset
{\Bbb L_\ep^\beta}$. This means the trajectories starting from $B$ at time 0 is always
attracted by the moving compact set ${\cal A_\ep }(\th_t\om)$.

\section{Impact of stochastic dynamic boundary conditions}

To investigate the impact of the stochastic dynamic boundary
conditions on overall evolution,  we estimate the Hausdorff
dimension of the random attractor associated with the stochastic
Cahn-Hilliard equations (\ref{1-5}) or (\ref{3-1}). We hope to
find a upper bound of the Hausdorff dimension of the random
attractor ${\cal A_\ep(\om)}$ by using Theorem \ref{dimension}.

Define a random map $S$ in ${\cal L}_\ep^\beta$ by
$$S(\om)=S_\ep(1,0;\om)=\ph(1,\om) ,$$
and an ergodic transformation $\th=\th_1 $. Then the random
attractor ${\cal A_\ep}(\om),\om\in\Om $ is a compact measurable
set invariant by $S$. In order to apply Theorem \ref{dimension},
we need to check that the three assumptions of this theorem hold.
To establish  the uniform differentiability of the random map $S$,
we first give some estimates.

Let $(\phi_i(t), \psi_i(t))(i=1,2)$ be two solutions
of the problem (\ref{3-1}) with the initial values $(\phi_i(0), \psi_i(0))=(\phi_i^0, \psi_i^0)\in \Bbb L_\ep^\beta $
and denote $(g(t), h(t))=(\phi_1(t)-\phi_2(t), \psi_1(t)-\psi_2(t))$. Then $(g(t), h(t))$ solves
\be\label{5-1}\ba{l}
(\ep+A^{-1})\pp_t g(t)=\De g-<\pp_\nu g>_\Ga-(\bar f(\phi_1)-\bar f(\phi_2)), x\in G,\\
g(0)=\phi_1^0-\phi^0_2,\\
\pp_t h(t)=\De_\|h-\la h-\pp_\nu g, x\in \Ga,\\
h(0)=\psi_1^0-\psi_2^0,\\
g|_\Ga=h. \ea\ee Taking the scalar product in $H$ of the first
equation of (\ref{5-1}) with $g$, by $m(g)=0$ and $f'(u)>-C$, we
get
\be\label{5-2} \pp_t|(g,h)|_{\Bbb
L_\ep^\beta}^2+2|(g,h)|_{\Bbb V}^2\le C|g|^2\le \frac12
|g|_1^2+C|g|_{-1}^2 , \ee
which leads to
\be\label{5-3}
|(g(t),h(t))|_{{\Bbb L}_\ep^\beta}^2\le e^{Ct}|(g(0),h(0))|_{{\Bbb
L}_\ep^\beta}^2, {\rm and } \int_0^1|(g(t),h(t))|_{\Bbb V}^2 dt\le
C|(g(0),h(0))|_{{\Bbb L}_\ep^\beta}^2.
\ee
The uniform differentiability of the random map $S$ can be stated as follows.
\begin{lea}\label{tdiff}
$S(\om)$ is almost surely uniformly differentiable on ${\cal A_\ep}(\om) $:
for $(\phi_0, \psi_0), (\phi_0+h, \psi_0+l) \in {\cal A_\ep}(\om)$, there exists $DS(\om; \phi_0,\psi_0)\in
{{\cal L} (\Bbb L^0_\ep)} $ such that
$$|S(\om)(\phi_0+h, \psi_0+l)-S(\om)(\phi_0, \psi_0)-DS(\om; \phi_0, \psi_0)(h, l)|_{\Bbb L_\ep^\beta}
\le \bar k(\om)|(h, l)|_{\Bbb L_\ep^\beta}^{1+\mu},$$ holds
P-a.s., where $\mu>0, \bar k(\om)\ge 1, E(\log \bar k(\om))<\iy$,
$DS(\om; \phi_0,\psi_0)(h, l)=(\Phi(1), \Psi(1))$, and $(\Phi(t),
\Psi(t))$ is a solution to the first variation system of the
system (\ref{3-1})
\be\label{5-4}\ba{l}
(\ep+A^{-1})\pp_t\Phi=\De \Phi-<\pp_\nu \Phi>_\Ga-\overline{f'(\phi)\Phi}, x\in G,\\
\Phi(0)=h,\\
\pp_t\Psi=\De_\|\Psi-\la\Psi-\pp_\nu\Phi, x\in \Ga,\\
\Psi(0)=l,\\
\Phi|_\Ga=\Psi,
\ea\ee
with $\phi(t)=S_\ep(t,0,;\om)\phi_0$ and $m(h)=0$.
\end{lea}

\noindent{\it Proof.} We write $(r(t), q(t))=(\phi_1(t)-\phi_2(t)-\Phi(t), \psi_1(t)-\psi_2(t)-\Psi(t))$,
where $(\phi_i(t), \psi_i(t))(i=1,2)$ are two solutions
to the problem (\ref{3-1}) with $(\phi_i(0), \psi_i(0))=(\phi_i^0, \psi_i^0)\in \cal A_\ep(\om)$ and $(\Phi(t), \Psi(t))$
satisfies the linear system (\ref{5-4}) with $\phi=\phi_2$ and
$(h, l)=(\phi_1^0-\phi_2^0, \psi_1^0-\psi_2^0)$. Then $(r(t), q(t))$ solves
\be\label{5-5}\ba{l}
(\ep+A^{-1})\pp_t r=\De r-<\pp_\nu r>_\Ga-\overline {f(\phi_1)}+\overline{f(\phi_2)}+\overline{f'(\phi_2)\Phi}, x\in G,\\
r(0)=0,\\
\pp_t q=\De_\|q-\la q-\pp_\nu r, x\in \Ga,\\
q(0)=0,\\
r|_\Ga=q. \ea\ee It is easy to obtain $m(\Phi)=m(r)=0$ since
$S_\ep(t)$ maps $\Bbb L_\ep^\beta$ into itself for $\fo \beta\in
\Bbb R$. Taking the scalar product in $H$ of the first equation
with $r$, we get \be\label{5-6} \pp_t |(r(t), q(t))|^2_{\Bbb
L_\ep^\beta}+2|(r(t), q(t))|^2_{\Bbb V}+2\la |q|_\Ga^2=
2(-\overline
{f(\phi_1)}+\overline{f(\phi_2)}+\overline{f'(\phi_2)\Phi}, r).\
\ee
By the H\"older inequality and the Sobolev embedding theorem,
for $s\in [\frac 65, 2]$ and its conjugate exponent $s^*\in
[2,6]$, the right hand side of the above inequality may be
estimated as
 \be\label{5-7}\ba{l}
\quad 2(-\overline {f(\phi_1)}+\overline{f(\phi_2)}+\overline{f'(\phi_2)\Phi}, r)\\
=2(-{f(\phi_1)}+ {f(\phi_2)}+ {f'(\phi_2)}(\phi_1-\phi_2-r), r)\\
\le C|r|^2+2|{f(\phi_1)}-{f(\phi_2)}-{f'(\phi_2)}(\phi_1-\phi_2)|_{L^s}|r|_{L^{s^*}}\\
\le |r|_1^2+C|r|_{-1}^2+C_1|{f(\phi_1)}-{f(\phi_2)}-{f'(\phi_2)}(\phi_1-\phi_2)|_{L^s}^2,
\ea\ee
where the fact $f'(\phi_2)\ge -C$ is used. We continue to estimate the third term of
(\ref{5-7}). For $0<\de<\frac 2s-1$, we have
$$\ba{l}
\quad |f(\phi_1)-f(\phi_2)-f'(\phi_2)(\phi_1-\phi_2)|_{L^s}^s\\
\le C|(|\phi_1|^{2p-3}+|\phi_2|^{2p-3}+1)(\phi_1-\phi_2)^2|^s_{L^s}\\
\le C\int_G(|\phi_1|+|\phi_2|+1)^{s(2p-2-\de)}|\phi_1-\phi_2|^{s(1+\de)} dx\\
\le C(|\phi_1|_{L^{s_1}}^{s(2p-2-\de)}+|\phi_2|_{L^{s_1}}^{s(2p-2-\de)}+1)|\phi_1-\phi_2|^{s(1+\de)},
\ea
$$
where $s_1=\frac{2s(2p-2-\de)}{2-s(1+\de)}$. Therefore, Eq.(\ref{5-6}) gives
\be\label{5-8}\ba{l}
\quad \pp_t |(r(t), q(t))|^2_{\Bbb L_\ep^\beta}\\
\le C|r|_{-1}^2+
C_1(|\phi_1|_{L^{s_1}}^{2(2p-2-\de)}+|\phi_2|_{L^{s_1}}^{2(2p-2-\de)}+1)|\phi_1-\phi_2|^{2(1+\de)}\\
\le C|r|_{-1}^2+
C_1(|\phi_1|_{L^{s_1}}^{2(2p-2-\de)}+|\phi_2|_{L^{s_1}}^{2(2p-2-\de)}+1)|\phi_1-\phi_2|_1^{1+\de}
|\phi_1-\phi_2|_{-1}^{1+\de}.
\ea\ee
By (\ref{5-3}) and (\ref{5-8}), we have
\be\label{5-9}\ba{rl}
|(r(1), q(1))|^2_{\Bbb L_\ep^\beta}&\le
C\int_0^1 (|\phi_1|_{L^{s_1}}^{2(2p-2-\de)}+|\phi_2|_{L^{s_1}}^{2(2p-2-\de)}+1)|\phi_1-\phi_2|_1^{1+\de}
ds |(h,l)|_{\Bbb L_\ep^\beta}^{1+\de }\\
&\le C\{\int_0^1 (|\phi_1|_{L^{s_1}}^{s_2}+|\phi_2|_{L^{s_1}}^{s_2}+1)ds |(h,l)|_{\Bbb L_\ep^\beta}^{2(1+\de)}\\
&:=k^2|(h,l)|_{\Bbb L_\ep^\beta}^{2(1+\de)},
\ea\ee
where $s_2=4(2p-2-\de)/(1-\de)$. Choosing $\bar k =\max\{k,1\}$, which satisfies
$E(\log\bar k(\om))<\iy$ by (\ref{4-11}), we conclude the proof of Lemma
\ref{tdiff}.
$\blacksquare$

\medskip

We continue to verify that assumptions (2) and (3) of Theorem \ref{dimension} hold in terms of
the following lemma.

\begin{lea}
For any $(\phi_0, \psi_0)\in {\cal A}(\om)$, there exist a random variable $\bar \al_1(\om)$ such that
\be\label{5-10}
\al_1(DS(\om; \phi_0,\psi_0))\le \bar \al_1(\om),\bar \al_1(\om)>1, E(\log \bar \al_1(\om))<\iy,
\ee
and another random variable $\bar\om_d(\om)$ such that
\be\label{5-11}
\om_d(DS(\om; \phi_0,\psi_0))\le \bar\om_d(\om), E\log \bar\om_d(\om)<0,
\ee
where the notations $\om_d$ and $\al_1$ are defined in Theorem \ref{dimension}.
\end{lea}

\noindent{\it Proof.} Multiplying the first equation of the variation system (\ref{5-1}) scalarly
in $H_0$ by $\Phi$, for $h\in H_0$, we have
\be\label{5-12}
\pp_t|(\Phi,\Psi)|_{\Bbb L_\ep^\beta}^2+|(\Phi, \Psi)|_{\Bbb V}^2+2\la|\Psi|_\Ga^2\le 2C|\Phi|_{-1}^2.
\ee
Hence, we get $|(\Phi(t),\Psi(t))|_{\Bbb L_\ep^\beta}\le
|(h,l)|_{\Bbb L_\ep^\beta}e^{Ct}$.
Since $\al_1(DS(\om; \phi_0, \psi_0))$ is equal to
the norm of $DS(\om; \phi_0, \psi_0)\in {{\cal L}({\Bbb L}^0_\ep)}$, we can choose a
constant $\bar \al_1(\om)=e^{C}$ such that
$$\al_1(DS(\om,\phi_0,\psi_0))\le \bar \al_1(\om), E(\log\bar\al_1(\om))<\iy, $$
which gives (\ref{5-10}). As for (\ref{5-11}), we first write
$$DS(\om; \phi_0, \psi_0)=\exp\{\int_0^1L(s,\phi(s),\psi(s))ds\},$$
and $L(s,\phi(s),\psi(s))(\Phi,\Psi)=((\ep+A^{-1})^{-1}(\De \Phi-<\pp_\nu\Phi>_\Ga-\overline{f'(\phi)\Phi}),
\De_\|\Psi-\la\Psi)$.
Following Ref. \cite{tem}, we have
$$\om_d(DS(\om; \phi_0,\psi_0))\le\sup_{{\xi_i\in H }\atop{|\xi_i|\le
1,i=1,\cdots,d}}\exp\{\int_0^1{\rm tr}( L(s,\phi(s),\psi(s))\circ
Q_d(s))ds\},$$ where $Q_d(s)$ is the orthonormal projector in $H_0\times H(\Ga)$
onto the space spanned by $(\Phi_1(s), \Psi_1(s))$,
$\cdots, (\Phi_d(s), \Psi_d(s))$, and $(\Phi_i(s), \Psi_i(s))$ is
the solution to the system (\ref{5-4})  with $(\Phi_i(0),\Psi_i(0)=(\xi_i, \eta_i)$.

Let $(\ph_i(s),\psi_i(s))\in {\Bbb V}, i\in {\Bbb N}$ be an
orthonormal basis of $H_0\times H(\Ga)$ such that
$Q_d(s)(H_0\times H(\Ga))={\rm
Span}[(\ph_1(s),\psi_1(s)),\cdots,(\ph_d(s),\psi_d(s))]$. Then
\be\label{5-13}\ba{l}
\quad {\rm tr}(L(s,\phi(s),\psi(s))\circ Q_d(s))\\
=\sum_{i=1}^d(L(s,u(s))(\ph_i(s),\psi_i(s),(\ph_i(s),\psi_i(s)))\\
=\sum_{i=1}^d (\De \ph_i-<\pp_\nu \ph_i>_\Ga-\overline{f'(\phi)\ph_i}, (\ep+A^{-1})^{-1}\ph_i)+
\sum_{i=1}^d(\De_\|\psi_i-\la\psi_i-\pp_\nu\ph_i,\psi_i)_\Ga\\
\le-C\sum_{i=1}^d|\ph_i|_2^2+\sum_{i=1}^d(-\overline{f'(\phi)\ph_i}, (\ep+A^{-1})^{-1}\ph_i)
-\sum_{i=1}^d|\psi_i|_{2,\Ga}^2-(\pp_\nu\ph_i,\psi_i)_\Ga\\
\le-C\sum_{i=1}^d|\ph_i|_2^2+\sum_{i=1}^d(-\overline{f'(\phi)\ph_i}, (\ep+A^{-1})^{-1}\ph_i)
-\sum_{i=1}^d|\psi_i|_{2,\Ga}^2+\sum_{i=1}^d|\psi_i|_{\frac32, \Ga}|\ph_i|\\
\le -C\sum_{i=1}^d|\ph_i|_2^2+\sum_{i=1}^d(-\overline{f'(\phi)\ph_i}, (\ep+A^{-1})^{-1}\ph_i)
+C_1d^2.
\ea\ee
Writing $ \zeta=\sum_{i=1}^d |\ph_i |^2$ and choosing H\"older exponents $p_n=1+4/n,
q_n=2+n/2,n<4$, by the generalized Sobolev-Lieb-Thirring inequality (Ref. \cite{tem})
$$|\zeta|_{L^{p_n}}^{p_n}\le C\sum_{i=1}^d |\ph_i |_2^2,\quad
\sum_{i=1}^d |\ph_i |_2^2\ge C_1d^{1+4/n}-C_2 d,$$
we have
$$\ba{l}
\quad {\rm tr}(L(s,\phi(s),\psi(s))\circ Q_d(s))\\
\le -C\sum_{i=1}^d |\ph_i |_2^2+ (|f'(\phi)|
\zeta^{1/2},(\sum_{i=1}^d |(\ep+A^{-1})^{-1}\ph_i |^2)^{1/2})+C_1d^2\\
\le -C\sum_{i=1}^d |\ph_i |_2^2+ |f'(\phi)|_{L^{q_n}}|\zeta|_{L^{p_n}}^{1/2}(\sum_{i=1}^d |\ph_i |_2^2)^{1/2}
+C_1d^2\\
\le -\frac C2 \sum_{i=1}^d |\ph_i |_2^2+C_1|f'(\phi)|_{L^{q_n}}^{q_n}+C_2d^2\\
\le -C_1 d^{1+4/n}+C_2 d++C_3d^2+ C|f'(\phi)|_{L^{q_n}}^{q_n}\\
\le -C_1 d^{1+4/n}+ C_2+C|f'(\phi)|_{L^{q_n}}^{q_n}. \ea$$ In the
end, we deduce \be\label{5-14} \om_d(DS(\om;
\phi_0,\psi_0))\le\exp\{ -C_1 d^{1+4/n}+
C_2+C\sup_{(\phi_0,\psi_0)\in \cal A_\ep(\om)} \int_0^1
|f'(\phi(t))|_{L^{q_n}}^{q_n} dt\}. \ee By the invariance of $\cal
A_\ep$ and (\ref{4-11}), the third term in the exponent of the
right hand side is finite.
Denoting by $\bar\om_d(\om)$  the right
hand side of (\ref{5-14}) and choosing $d$ such that
\be\label{5-15} -C_1 d^{1+4/n}+ C_2+C\sup_{(\phi_0,\psi_0)\in \cal
A_\ep(\om)}\int_0^1 |f'(\phi(t))|_{L^{q_n}}^{q_n} dt< 0, \ee then
we have $\om_d(DS)\le \bar\om_d(\om)$ and $E(\log(\bar\om_d))<0$.

In conclusion, we have the following result.
\begin{thm}
If there exists $d$ such that (\ref{5-15}) holds, then Hausdorff
dimension of the random attractor for the   Cahn-Hilliard system
under stochastic dynamic boundary condition (\ref{1-5}) is finite
and the dimension is bounded by $d$, i.e.,    $d_H({\cal
A_\ep(\om)})<d, \; P-a.s.$
\end{thm}

\begin{rek} (\textbf{Impact of stochastic dynamic boundary conditions}.)
A interesting problem is that if the dynamic  boundary condition
in Eq.(\ref{1-5}) is replaced by \be\label{5-16}
\frac1{\ep_0}d\psi=(\De_\|\psi-\la\psi-\pp_\nu\phi)dt+\si_2
dW^{(2)}, \,  x\in \Ga , \ee where $\ep_0>0$ is a scaling
parameter modeling the ``intensity" of the time-derivative
component in the dynamic boundary condition,  then (\ref{5-15})
has to be rewritten as
$$
-C_1 d^{1+4/n}+ C_2\ep_0^{\frac{4+n}{4-n}}+C_3
+C\sup_{(\phi_0,\psi_0)\in \cal A_\ep(\om)}\int_0^1 |f'(\phi(t))|_{L^{q_n}}^{q_n} dt< 0,
$$
which implies $d \sim \ep_0^{\frac n{4-n}}$ and describes the
impact of the dynamical boundary conditions on the global
dynamics.

We may think $\frac1{\ep_0} $ as the intensity (relative
importance) of the time derivative component in the dynamic
boundary conditions (\ref{5-16}). If $\ep_0$ is very big (i.e.,
the time derivative term in the stochastic  dynamic boundary
conditions is very small), then the dimension is very big (Note
that $ n=1, 2, 3$). This    says that the dimension of the random
attractor increases as $\frac1{\ep_0} $ decreases, namely, as the
time derivative term   becomes weaker.

However,   in the limiting case $\ep_0=\iy$, i.e., the stochastic
dynamic boundary conditions reduce to the stochastic static
boundary conditions, \be \label{5-17}
0=(\De_\|\psi-\la\psi-\pp_\nu\phi)dt+\si_2 dW^{(2)}, \,  x\in \Ga
, \ee the dimension \emph{does not } tend to infinity. Instead, we
see that the stochastic \emph{static} boundary conditions
(\ref{5-17}) do not have impact on the dimension. In fact,
 we can remove
$\sum_{i=1}^d(\De_\|\psi_i-\la\psi_i-\pp_\nu\ph_i,\psi_i)_\Ga$
from the crucial equation (\ref{5-13}), and thus, the stochastic
static boundary conditions do not affect the later derivations and
thus have no    impact on the Hausdorff dimension of the random
attractor.
\end{rek}


\end{document}